\documentclass[12pt,oneside]{amsart}

\usepackage{amscd,amsmath,amsfonts,amssymb,amsthm,array}
\usepackage[all]{xy}

\usepackage[a4paper]{geometry}

\newcommand{\C}{\ensuremath{\mathbb C}}
\newcommand{\Ha}{\ensuremath{\mathbb H}}

\newcommand{\R}{\ensuremath{\mathbb R}}
\newcommand{\Q}{\ensuremath{\mathbb Q}}
\newcommand{\Z}{\ensuremath{\mathbb Z}}

\newcommand{\CP}[1]{\ensuremath{{\C\mathbb P}^{#1}}}

\newcommand{\integer}[1]{\ensuremath{{\mathcal O}_{#1}}}
\newcommand{\unit}[1]{\ensuremath{{\mathcal O}_{#1}^*}}
\newcommand{\unitpos}[1]{\ensuremath{{\mathcal O}_{#1}^{*,+}}}
\newcommand{\norm}[1]{\ensuremath{\text{\upshape \rmfamily Nm}(#1)}}
\newcommand{\normover}[3]{\ensuremath{\text{\upshape \rmfamily Nm}_{#1/#2}(#3)}}
\newcommand{\coker}{\operatorname{coker}}
\newcommand{\Gal}{\operatorname{Gal}}

\usepackage{graphicx}
\newcommand{\turndown}[1]{%
  \rotatebox[origin=c]{270}{\ensuremath#1}}
\newcommand{\hookdownarrow}{\turndown{\hookrightarrow}}
\newcommand{\twoheaddownarrow}{\turndown{\twoheadrightarrow}}

\newcommand{\Aut}[1]{\ensuremath{\text{\upshape \rmfamily Aut}(#1)}}

\newcommand{\Hmt}[1]{\ensuremath{\text{\upshape \rmfamily Hmt}(#1)}}
\newcommand{\Hom}[2]{\ensuremath{\text{\upshape \rmfamily Hom}(#1,#2)}}

\newcommand{\lckrk}[1]{\ensuremath{\text{\upshape \rmfamily r}_{#1}}}
\newcommand{\rk}[1]{\ensuremath{\text{\upshape \rmfamily Rank}(#1)}}

\newcommand{\iso}{\ensuremath{\simeq}}
\newcommand{\lck}{\textsc{lck}}

\newcommand{\lee}{\ensuremath{\omega}}

\newcommand{\ug}{\stackrel{\rm def}{=}}
\newcommand{\fund}{\ensuremath{\Omega}}
\newcommand{\map}[3]{\mbox{${#1}\colon{#2}\to{#3}$}}

\newcommand{\st}{\ensuremath{\text{ such that }}}
\newcommand{\dismap}[5]{
\begin{gather*}
#2 \longrightarrow #3\\
#1(#4)\ug #5
\end{gather*}
}

\renewcommand{\Im}[1]{\ensuremath{\text{\upshape \rmfamily Im}#1}}

\numberwithin{equation}{section}

\newtheorem{te}{Theorem}[section]
\newtheorem{lm}[te]{Lemma}

\newtheorem*{te*}{Theorem}

\theoremstyle{definition}

\newtheorem{de}[te]{Definition}
\newtheorem{re}[te]{Remark}
\newtheorem{ex}[te]{Example}

\newcommand{\EndDim}{\ensuremath{\nopagebreak\hfill\blacksquare}}

\def\ra{\rightarrow}

\newenvironment{D}[1][]{{\nopagebreak\noindent\em Proof#1: }}{\EndDim}

\begin{document}

\title[Examples of non-trivial rank in LCK geometry]{Examples of non-trivial rank in locally conformal K\"ahler geometry}

\author{Maurizio Parton, Victor Vuletescu}

\address{Maurizio Parton,\newline Universit\`a di Chieti--Pescara\\
Dipartimento di Scienze, viale Pindaro $87$, I-65127 Pescara, Italy.}
\email{parton@sci.unich.it}
\address{Victor Vuletescu,\newline University of Bucharest, Faculty of Mathematics and Informatics,
Academiei st. 14,
70109 Bucharest, Romania.}
\email{vuli@fmi.unibuc.ro}

\begin{abstract}
We consider locally conformal K\"ahler geometry as an equivariant, homothetic K\"ahler geometry $(K,\Gamma)$. We show that the de Rham class of the Lee form can be naturally identified with the homomorphism projecting $\Gamma$ to its dilation factors, thus completing the description of locally conformal K\"ahler geometry in this equivariant setting. The rank $\lckrk{M}$ of a locally conformal K\"ahler manifold is the rank of the image of this homomorphism. Using algebraic number theory, we show that $\lckrk{M}$ is non-trivial, providing explicit examples of locally conformal K\"ahler manifolds with $1\nless\lckrk{M}\nless b_1$. As far as we know, these are the first examples of this kind. Moreover, we prove that locally conformal K\"ahler Oeljeklaus-Toma manifolds have either $\lckrk{M}=b_1$ or $\lckrk{M}=b_1/2$.
\end{abstract}

\maketitle


\section{Introduction}

For many reasons, K\"ahler manifolds are considered the most interesting objects of complex geometry. However, strong topological properties -like formality- even Betti numbers of odd index and others, obstruct the existence of K\"ahler metrics on many compact manifolds, some of them very simple ones, like the Hopf or Kodaira surfaces. From the Riemannian viewpoint, the natural place to look for metrics with a given property is a conformal class. When this is not possible, then local metrics with the said property can be searched for, subject to some condition on the overlaps.

This is exactly the way Izu Vaisman arrived to the notion of \emph{locally conformal K\"ahler} (briefly, \lck) metric \cite{VaiLCA}. The original definition puts the accent on a fixed metric which is locally conformal with local K\"ahler ones. Equivalently, it requires the existence of a closed one-form (the Lee form) which, together with the fundamental two-form, generates a differential ideal. On the other hand, any metric globally conformal with a \lck\ metric is again \lck. This allows talking about a \lck\ structure, in which no metric is fixed and only the cohomology class of the Lee form is given. This understanding of \lck\ geometry is consistent with the fact that any K\"ahler cover of a \lck\ manifold bears a K\"ahler metric with respect to which the covering group acts by holomorphic homotheties. \lck\ geometry can thus be seen as the pair $(K,\Gamma)$ of a K\"ahler manifold and a group of holomorphic homotheties. This viewpoint has been suggested in \cite{GOPLCK}, and then developed in \cite{GOPRVS}, where two key notions were introduced: the \emph{presentation} (in this paper called \lck-presentation), which is the pair described above, and the \emph{rank} of the subgroup of $\R^+$ given by dilation factors of $\Gamma$, which measures the ``true'' homothety part of the group.

In the present paper we go a bit further, showing that the Lee form can also be read in these terms. This completes the description of \lck\ geometry in terms of presentations. Moreover, we show that the examples of \lck\ manifolds constructed in \cite{OeTNKC} have highly non-trivial rank: their rank is either equal to the first Betti number or to half of it. In particular, this provides a first example of \lck\ manifold of rank $\neq 1$ and strictly less than $b_1$.

The structure of the paper is as follows. In Section $2$, we recall the basic definitions and properties of \lck\ geometry, presentations and rank. In Section $3$ we show how the Lee form can be reconstructed from a presentation. Section $4$ is devoted to a detailed description of the complex manifolds defined by Oeljeklaus-Toma in \cite{OeTNKC}, and to the computation of their first Betti number. In Section $5$, we recall how Oeljeklaus-Toma manifolds can be \lck-presented, in terms of a global potential, and we compute the dilation factors of $\Gamma$. Then we prove the following Theorem:
\begin{te*}
Let $M$ be an \lck\ Oeljeklaus-Toma manifold. Then its rank is either $b_1(M)$ or $b_1(M)/2$.
\end{te*}
Using this Theorem we then compute explicit examples of \lck\ manifolds with non-trivial rank.

Since Section $4$ and $5$ makes strong use of tools from Algebraic Number Theory, in Section $6$ we make a short summary of these tools.

\section{\textsc{lck}-presentations for complex manifolds}

For convenience of the reader, we here briefly review notation established in \cite{GOPRVS}.

Let $M$ be a complex manifold. A \emph{locally conformal K\"ahler} metric is a conformal class $[g]$ of Hermitian metrics on $M$ such that $[g]$ is given locally by K\"ahler metrics. The conformal class $[g]$ corresponds to a unique de Rham cohomology class $[\lee_g]\in H^1(M)$, whose representative 
$\lee_g$ is defined as the unique closed 1-form satisfying $d\fund_g=\lee_g\wedge\fund_g$, where $\fund_g$ denotes the fundamental form of $g$. The 1-form $\lee_g$ is called the \emph{Lee form} of $g$.

Taking into account that a locally conformal K\"ahler metric on a manifold of K\"ahler type must be globally conformal K\"ahler \cite{VaiLCK}, it is a trivial task to show that a complex manifold $M$ (of complex dimension at least $2$) admits a locally conformal K\"ahler metric if and only if there is a complex covering space $K$ of K\"ahler type such that $\pi_1(M)$ acts on $K$ by holomorphic homotheties with respect to the K\"ahler metric.

More explicitly, if $g$ is a locally conformal K\"ahler metric on $M$, and $\lee_g$ its Lee form, then the pull-back of $\lee_g$ to any K\"ahler covering $K$ of $M$ is exact, say $\tilde{\lee}_g=df$. Denoting by $\tilde{g}$ a lift of $g$ to $K$, then $e^{-f}\tilde{g}$ turns out to be a K\"ahler metric on $K$ such that $\pi_1(M)$ acts on it by holomorphic homotheties.
According to the fact that the K\"ahler metric $e^{-f}\tilde{g}$ is defined up to homotheties (because $f$ is defined up to a constant), usually in locally conformal K\"ahler geometry one is interested in the \emph{homothety class} of a K\"ahler manifold.

The above discussion motivates the following definitions, first given in \cite{GOPRVS}. For the notion of minimal cover in the more general setting of conformal geometry, see also \cite{BeMHTW}.

\begin{de}
Let $K$ be a homothetic K\"ahler manifold and $\Gamma$ a discrete Lie group of biholomorphic homotheties acting freely and properly discontinously on $K$.
\begin{itemize}
\item The pair $(K,\Gamma)$ is called a \emph{\textsc{lck}-presentation}.
\item If $M$ is a complex manifold and $M=K/\Gamma$ as complex manifolds, $(K,\Gamma)$ is called a \emph{\textsc{lck}-presentation for} $M$.
\item If $\Gamma$ does not contain isometries other than the identity, then $(K,\Gamma)$ is called \emph{minimal}, and if $K$ is simply connected then $(K,\Gamma)$ is called \emph{maximal}.
\end{itemize}
\end{de}

\begin{re}
Given a complex manifold $M$, the statement ``$(K,\Gamma)$ is a \textsc{lck}-presentation for $M$'' is just a shortcut for ``$K$ is a complex covering space of $M$, and $\Gamma$ are its covering transformations, and there is a K\"ahler metric on $K$ which is conformally equivalent to a $\Gamma$-invariant metric''. Due to the $1$-$1$ correspondence existing between locally conformal K\"ahler manifolds and minimal presentations, we will often abuse of this language by saying ``the locally conformal K\"ahler manifold $(K,\Gamma)$''.
\end{re}

In a homothetic K\"ahler manifold $K$ we denote by $\Hmt{K}$ the group of its biholomorphic homotheties, and by
\[
\map{\rho_K}{\Hmt{K}}{\R^+}
\]
the group homomorphism associating to a homothety its dilation factor.
For any locally conformal K\"ahler manifold $M$, \textsc{lck}-presented as $(K,\Gamma)$, the rank of the free abelian group $\rho_K(\Gamma)$ depends only on $M$ \cite[Proposition~2.10]{GOPRVS}.

\begin{de}
The rank of $\rho_K(\Gamma)$ is called \emph{the rank of $M$}, and is denoted by \lckrk{M}.
\end{de}

\begin{re}
The rank $\lckrk{M}$ measures ``how much'' the locally conformal K\"ahler manifold is far from the K\"ahler geometry.
\end{re}

\section{The Lee form}

Let $M$ be a locally conformal K\"ahler manifold \textsc{lck}-presented as $(K,\Gamma)$. The question if the de Rham class of any Lee form of $M$ can be completely described in terms of \textsc{lck}-presentations has been left open in \cite{GOPRVS}. In this Section we fill this gap.

Consider the following exact sequence:
\[
1 \rightarrow \pi_1(K) \rightarrow \pi_1(M)\rightarrow \Gamma \rightarrow 1
\]
Recalling that the abelianization doesn't preserve exactness on the left, we get:
\begin{equation}\label{H1}
H_1(K,\Z) \rightarrow H_1(M,\Z)\rightarrow \frac{\Gamma}{[\Gamma,\Gamma]} \rightarrow 0
\end{equation}
Using the universal coefficient theorem for cohomology and the de Rham theorem, we can translate the above sequence in de Rham cohomology language:
\begin{equation}\label{i}
0 \rightarrow \Hom{\frac{\Gamma}{[\Gamma,\Gamma]}}{\R} \stackrel{i}{\rightarrow} H^1_{dR}(M) \rightarrow H^1_{dR}(K)
\end{equation}

Now, observe that $\Gamma\subset\Hmt{K}$ (by definition), and $[\Gamma,\Gamma]\subset\ker{\rho_K}$ (because $\R$ is abelian). We thus obtain a homomorphism from $\frac{\Gamma}{[\Gamma,\Gamma]}$ to $\R$ by
\begin{equation}\label{rho_K}
\frac{\Gamma}{[\Gamma,\Gamma]}\stackrel{\frac{\rho_K|\Gamma}{[\Gamma,\Gamma]}}{\longrightarrow}\R^+\stackrel{\log}{\longrightarrow}\R
\end{equation}

For the sake of simplicity, we still denote by $\rho_K$ this element of $\Hom{\frac{\Gamma}{[\Gamma,\Gamma]}}{\R}$. We are now ready to state the following Theorem.

\begin{te}\label{lee}
Let $M$ be a locally conformal K\"ahler manifold \textsc{lck}-presented as $(K,\Gamma)$, and let $[\lee]\in H^1_{dR}(M)$ be its Lee form. Let $i$ be the map given by \eqref{i}, and $\rho_K$ the element of $\Hom{\frac{\Gamma}{[\Gamma,\Gamma]}}{\R}$ given by \eqref{rho_K}. Then
\[
[\lee]=i(\rho_K)
\]
\end{te}

\begin{D}
Denote by $p$ the projection from $K$ to $M$, and by $g$ the Riemannian metric on $M$ associated to $\lee$. Thus, $p^*g$ is a $\Gamma$-invariant metric on $K$, $p^*\lee=df$ is an exact 1-form on $K$, and the metric $g_K=e^{-f}p^*g$ on $K$ is K\"ahler.

For any $\gamma\in\Gamma$, denote by $[\gamma]$ the corresponding element of $\frac{\Gamma}{[\Gamma,\Gamma]}$. We then have:
\[
\gamma^*g_K=\gamma^* e^{-f}p^*g=e^{-f\circ\gamma}\gamma^*p^*g=e^{-f\circ\gamma}p^*g=e^{-f\circ\gamma+f}e^{-f}p^*g=e^{-f\circ\gamma+f}g_K
\]
and thus (remember that $\rho_K$ is defined as in \eqref{rho_K}):
\begin{equation}\label{rhogamma}
\rho_K([\gamma])=-f\circ\gamma+f
\end{equation}
Remark in particular that since $\gamma\in\Hmt{K}$, then $-f\circ\gamma+f$ is constant.

To prove the claim, we thus need to show that $\lee([\alpha])=\rho_K(\gamma_\alpha^{-1})$ for every loop $\alpha$ in $M$. As for $\lee([\alpha])$, we have:
\[
\lee([\alpha])=\int_\alpha\lee=\int_{\tilde{\alpha}_{y_0}}p^*\lee=\int_{\tilde{\alpha}_{y_0}}df=f(\tilde{\alpha}_{y_0}(1))-f(\tilde{\alpha}_{y_0}(0))
\]
As for $\rho_K(\gamma_\alpha^{-1})$, we use \eqref{rhogamma}:
\[
\begin{split}
\rho_K(\gamma_\alpha^{-1})=-f\circ\gamma_\alpha^{-1}+f&=-f\circ\gamma_\alpha^{-1}(\tilde{\alpha}_{y_0}(1))+f(\tilde{\alpha}_{y_0}(1))\\
&=-f(\tilde{\alpha}_{y_0}(0))+f(\tilde{\alpha}_{y_0}(1))
\end{split}
\]
and this proves the claim.
\end{D}

\begin{re}
In the proof of Theorem~\ref{lee} we have shown that the rank can be defined in terms of the Lee form $\lee$, as the rank of the image of the natural map
\begin{gather*}
H_1(M,\Z) \longrightarrow \R\\
\alpha\longmapsto \int_\alpha \lee
\end{gather*}
\end{re}

\begin{re}
The rank $\lckrk{M}$ satisfies $0\le\lckrk{M}\le b_1(M)$, and $\lckrk{M}=0$ if and only if $M$ is globally conformal K\"ahler.
\end{re}

\section{Oeljeklaus-Toma manifolds}\label{secOT}

In their beautiful paper \cite{OeTNKC}, the authors construct locally conformal K\"ahler manifolds using tools from Algebraic Number Theory, which are summarized in Section~\ref{background}. In this section, we assume these tools are known.

We will denote by $F$ an algebraic number field, by \integer{F} the ring of algebraic integers of $F$ and by $\unit{F}$ the multiplicative group of units of $\integer{F}$. If $[F:\Q]=n=s+2t$ is the degree of $F$ over $\Q$, we denote by $\{\map{\sigma_i}{F}{\C}\}_{i=1,\dots,n}$ the complex embeddings of $F$, where the first $s$ embeddings are real, and the last $2t$ satisfy $\sigma_{s+i}=\bar\sigma_{s+i+t}$. The units which are positive in all real embeddings of $F$ are denoted by $\unitpos{F}$.

We are now ready to describe Oeljeklaus-Toma construction. For details, see \cite{OeTNKC}.

All together, the embeddings $\sigma_i$ give the natural map
\dismap{\sigma}{F}{\C^{s+t}}{x}{(\sigma_1(x),\dots,\sigma_{s+t}(x))}

The image $\sigma(\integer{F})$ is a lattice of rank $n$ in $\C^{s+t}$ \cite[Proposition 4.26]{MilANT}, and in this way we get a properly discontinuous action of $\integer{F}$ on $\C^{s+t}$ given by translations. We denote this action by $T$: if $a\in\integer{F}$ and $(z_1,\dots,z_{s+t})\in\C^{s+t}$, then
\begin{equation}\label{tra}
T_a(z_1,\dots,z_{s+t})\ug(z_1+\sigma_1(a),\dots,z_{s+t}+\sigma_{s+t}(a))
\end{equation}

We also have a multiplicative action of $\unit{F}$ on $\C^{s+t}$, denoted by $R$ (as in ``rotation''): if $u\in\unit{F}$, then
\begin{equation}\label{rot}
R_u(z_1,\dots,z_{s+t})\ug(\sigma_1(u)z_1,\dots,\sigma_{s+t}(u)z_{s+t})
\end{equation}

Pairs $(a,u)$ act then on $\C^{s+t}$ by $T_a\circ R_u$, and from $T_a\circ R_u\circ T_b\circ R_v=R_{uv}\circ T_{a+ub}$ one gets $(a,u)(b,v)=(a+ub,uv)$. In other words, the inclusion $\unit{F}\integer{F}\subset\integer{F}$ defines a semidirect product $\integer{F}\rtimes\unit{F}$ acting on $\C^{s+t}$.

Since for any $(a,u)\in\integer{F}\rtimes\unit{F}$ the equation $T_a(R_u(z_1,\dots,z_{s+t}))=(z_1,\dots,z_{s+t})$ has one solution $\sigma(\frac{a}{1-u})$, the action is not free, with fixed point set contained in $\sigma(F)\subset\R^s\times\C^t$. Thus, consider the upper complex half-plane $\Ha$ of complex numbers with strictly positive imaginary part, and observe that $\Ha^s$ is $\unitpos{F}$-invariant. Hence, $\integer{F}\rtimes\unitpos{F}$ acts freely on $\Ha^s\times\C^t$.

The above action is free, but not properly discontinuous in general (naively, one could say that ``there are too many generators'' for the group to act so: indeed $\integer{F}$ has rank $s+2t$, $\unitpos{F}$ has rank $s+t-1$ so there are $2s+3t-1$ generators, while we expect this number to be $2s+2t$). Still, Oeljeklaus and Toma show that one can always find a suitable subgroup $U\subset\unitpos{F}$ such that the action of $\integer{F}\rtimes U$ is properly discontinuous and moreover, the quotient is compact. Subgroups of these kind are called \emph{admissible}, and it is furthermore shown that when $t=1$ every subgroup of finite index of $\unitpos{F}$ is admissible.

Both $\integer{F}$ and $\unitpos{F}$ act holomorphically on $\Ha^s\times\C^t$, so the quotient inherits a complex structure.

\begin{de}
Given a finite field extension $F$ of $\Q$ and an admissible subgroup $U\subset\unitpos{F}$, the compact complex manifold
\[
M(F,U)\ug\frac{\Ha^s\times\C^t}{\integer{F}\rtimes U}
\]
is called an \emph{Oeljeklaus-Toma manifold}.
\end{de}

The first Betti number of $M(F,U)$ is computed in \cite{OeTNKC}. Since we also need this fact, we include here a different proof than the original one, which makes no use of group cohomology, spectral sequences and Hurewicz's Theorem.
\begin{te}\cite[Proposition 2.3]{OeTNKC}
Let $M=M(F,U)$ be an Oeljeklaus-Toma manifold. Then $b_1(M)=s$.
\end{te}

\begin{D}
We identify $\pi_1(M)$ with the deck transformation group $\integer{F}\rtimes U$, which is generated by $\{T_a,R_u\}_{a\in\integer{F},u\in U}$ see \eqref{tra}, \eqref{rot}. Since $\pi_1(M)/\integer{F}\iso U$ is abelian, and $H_1(M,\Z)$ is the maximal abelian quotient of $\pi_1(M)$, we have a commutative diagram:
\[
\begin{array}{ccccccccc}
1 & \ra & [\pi_1(M),\pi_1(M)] & \ra & \pi_1(M) & \ra & H_1(M,\Z) & \ra & 1\\
 & & i \hookdownarrow & & \parallel & & p\twoheaddownarrow\\
1 & \ra & \integer{F} & \ra & \pi_1(M) & \ra & U & \ra & 1\\
\end{array}
\]

From $p$ we get that $\rk{H_1(M,\Z)}=\rk{U}-\rk{\ker p}$. But the Snake Lemma gives $\ker p\iso \coker i$, thus it is enough to show that $\coker i=\integer{F}/[\pi_1(M),\pi_1(M)]$ is finite.

By direct computation we see $[T_a,R_u]=T_{(1-u)a}$, for any $a\in\integer{F}$ and any $u\in U$. In particular, this shows that for any $u\in U$, the principal ideal $(1-u)\integer{F}$ is a subgroup of $[\pi_1(M),\pi_1(M)]$. But if $u\neq 1$, then $(1-u)\integer{F}$ has finite index, as $\integer{F}$ is a Dedekind ring.
\end{D}

Since for the rest of this paper we will be concerned only with the case $t=1$ and $U=\unitpos{F}$, we skip the details about $U$.

\section{The rank of Oeljeklaus-Toma manifolds}\label{secpot}

The following result was the starting point for this paper.

\begin{te}\label{mainoe}\cite[Page 7]{OeTNKC}
Consider an Oeljeklaus-Toma manifold $M(F,\unitpos{F})$, with $t=1$ and $s>0$.
\begin{enumerate}
\item\label{pos} The real function
\begin{equation}\label{potential}
\phi(z)\ug\prod_{j=1}^s\frac{i}{z_j-\bar{z}_j}+\vert z_{s+1}\vert^2
\end{equation}
is a global K\"ahler potential on $\Ha^s\times\C$.
\item\label{equ} When $\Ha^s\times\C$ is equipped with the K\"ahler metric $i\partial\bar\partial \phi$ given above, the pair $(\Ha^s\times\C,\integer{F}\rtimes\unitpos{F})$ is a \textsc{lck}-presentation for $M(F,\unitpos{F})$.
\end{enumerate}
\end{te}

\begin{D}
To prove \eqref{pos}, one has to show that $\phi$ is strictly plurisubharmonic. By  direct calculation, we see that
\[
(\partial_{z_l}\partial_{\bar z_k}\phi)=
\left(
\begin{array}{cccc}
(\partial_{z_l}\partial_{\bar z_k}\phi_1)&0\\
0&2
\end{array}
\right)
\]
where
\[
\phi_1=\prod_{j=1}^s\frac{i}{z_j-\bar{z}_j}
\]
thus it suffices to look only at $(\partial_{z_l}\partial_{\bar z_k}\phi_1)$.

One has
\begin{equation}
\begin{gathered}\label{form}
\partial_{z_l}\partial_{\bar z_k}\phi_1=\frac{-1}{(z_k-\bar z_k)(z_l-\bar z_l)}\phi_1\qquad l\neq k=1,\dots,s\\
\partial_{z_k}\partial_{\bar z_k}\phi_1=\frac{-2}{(z_k-\bar z_k)^2}\phi_1\qquad k=1,\dots,s
\end{gathered}
\end{equation}
thus $(\partial_{z_l}\partial_{\bar z_k}\phi_1)$ is proportional to the matrix
\[
A=
\left(
\begin{array}{cccccc}
\frac{2}{4y_1^2} & \frac{1}{4y_1y_2} & \dots & \frac{1}{4y_1y_k}\\
\frac{1}{4y_2y_1} & \frac{2}{4y_2^2} & \dots & \frac{1}{4y_2y_k}\\
\cdot & \cdot & \dots & \cdot \\
\frac{1}{4y_ky_1} & \frac{1}{4y_ky_2} & \dots & \frac{2}{4y_k^2}
\end{array}
\right)
\]
But $A$ is the sum between a diagonal, positive definite matrix, and a positive semidefinite one, thus $A$ is positive definite.

Alternatively, one can directly notice that $\frac{i}{z-\bar{z}}$ defines a K\"ahler potential in $\Ha$ and then use \cite[Theorem~5.6]{DemCAD}.

To prove \eqref{equ}, one has to show that $(\integer{F}\rtimes\unitpos{F})$ act by homoteties on $\Ha^s\times\C$. Let $a\in\integer{F}$, $u\in\unitpos{F}$, and consider the generators $T_a$, $R_u$ of $(\integer{F}\rtimes\unitpos{F})$ given by \eqref{tra}, \eqref{rot}. Then, using \eqref{form} above and the fact that the embeddings $\{\sigma_j\}_{j=1,\dots,s}$ are real, one obtains $T_a^*(i\partial\bar\partial \phi)=i\partial\bar\partial \phi$, whereas using \eqref{prodsigma} one obtains
\[
R_u^*(\phi)=\vert\sigma_{s+1}(u)\vert^2 \phi
\]
that is, $R_u$ acts by homotheties on the potential itself.
\end{D}

\begin{re}
The K\"ahler potential $\phi$ given by \eqref{potential} corrects a minor typo present in the original paper. We acknowledge a useful email exchange with Matei Toma.
\end{re}

\begin{re}\label{image}
In the proof of Theorem~\ref{mainoe}, we have shown that the rank of $M(F,\unitpos{F})$ is the rank of the multiplicative subgroup of $\R^+$ given by
\[
\{ \vert \sigma_{s+1}(u) \vert^2 \st u\in\unitpos{F}\}
\]
(see also the proof of \cite[Proposition~2.9]{OeTNKC}).
\end{re}


The following result describes the rank of Oeljeklaus-Toma manifolds.
\begin{te}\label{propodd}
Let $F$ be a number field with $s>0$ real embeddings $\sigma_1,\dots,\sigma_s:F\rightarrow\R$ and exactly two non-real embeddings $\sigma_{s+1},\bar\sigma_{s+1}:F\rightarrow\C$. Let $M=M(F,\unitpos{F})$ be an Oeljeklaus-Toma manifold with the locally conformal K\"ahler structure given by Theorem \ref{mainoe}. Let $n=[F:\Q]$, so that $n=s+2$, $\dim_\C M=n-1$ and $b_1(M)=n-2$.
\begin{enumerate}
\item\label{odd} If $n$ is odd, then $M$ has maximal rank, that is to say, $\lckrk{M}=b_1(M)=n-2$.
\item\label{even} If $n$ is even, then either $M$ has again maximal rank or $\lckrk{M}=\frac{b_1(M)}{2};$ the last situation occurs if and only if $F$ is a quadratic extension of a totally real number field.
\end{enumerate}
\end{te}

\begin{D}
\eqref{odd} By Remark \ref{image}, it is enough to show that the map
\begin{equation}\label{hmtfactors}
u\mapsto\vert\sigma_{s+1}(u)\vert^2
\end{equation}
is injective.
Let $u\in\unitpos{F}$ be a unit with $\vert\sigma_{s+1}(u)\vert^2=1$ (that is to say, $\vert\sigma_{s+1}(u)\vert=1$), and let $P_u$ be its minimal polynomial over $\Q$. Since $p\ug\deg P_u=[\Q(u):\Q]$ divides $n=[F:\Q]$ and $n$ is odd, we see that also $p$ is odd. So $P_u$ is given by
\[
P_u(X)=X^p+a_{p-1}X^{p-1}+\dots+a_1X+a_0
\]
Moreover, using \eqref{char}, \eqref{norm} and \eqref{prodsigma} we obtain
\[
a_0^{n/p}=-\norm{u}=-1
\]
and since $a_0\in\Z$, we get $a_0=-1$.
Now observe that $\vert\sigma_{s+1}(u)\vert=1$ implies $\bar\sigma_{s+1}(u)=\frac{1}{\sigma_{s+1}(u)}$. But $\bar\sigma_{s+1}(u)$ is a root of $P_u$, hence $P_u\left(\frac{1}{\sigma_{s+1}(u)}\right)=0$. This means that $\sigma_{s+1}(u)$ satisfies the equation of degree $p$
\[
1+a_{p-1}X+\dots+a_1X^{p-1}-X^p=0
\]

Thus, the uniqueness of the minimal polynomial forces $a_k=-a_{p-k}$ for all $k$: but then $P_u(1)=0$, hence $u=1$.

\eqref{even} We consider the case when $M$ is not of maximal rank. This means that the map \eqref{hmtfactors} is not injective, so that there exists a unit $u\in\unitpos{F}$, $u\neq 1$, such that $\vert\sigma_{s+1}(u)\vert=1$. We claim that $\deg P_u$ is even. In fact, if $\deg P_u$ was odd, by \eqref{char} we would get
\[
a_0^{n/p}=\norm{u}=1
\]
that is, $a_0=\pm 1$. If $a_0=-1$, arguing the same as in point \eqref{odd} above, we would have $u=1$, whereas if $a_0=1$, we would have $P_u(-1)=0$, that is, $u=-1\not\in\unitpos{F}$, which is a contradiction in both cases. Thus, $\deg P_u$ is even.

\begin{lm}\label{primitive}
If $F$ admits exactly $2$ complex non-real embeddings, then any proper intermediate field extension $F\supsetneq E\supsetneq\Q$ is totally real.
\end{lm}

\begin{D}
Assume $E$ is not totally real, and let $\varsigma$ a complex non-real embedding. Let $d\ug [F:E]$. Then $F$ admits $d$ embeddings fixing $E$ pointwise, and their composition with $\varsigma$ gives $2d$ complex non-real embeddings of $F$. Thus $d=1$, and $F=E$.
\end{D}

In what follows, we need to assume that $u$ is non-real: since we can always replace $u$ with $\sigma_{s+1}(u)$, this assumption is not restrictive. Then, using Lemma~\ref{primitive}, we get $F=\Q(u)$.

Consider now the intermediate field extension $E\ug\Q(u+\frac{1}{u})$. Clearly, $F\supsetneq E\supsetneq\Q$. Using again Lemma~\ref{primitive}, we get that $E$ is totally real, whereas from $\left(u+\frac{1}{u}\right)u=u^2+1$ we get $[F=\Q(u):E]=2$.

It remains to check that the rank in this case is $\frac{b_1(M)}{2}$.
For any unit $u\in\unitpos{F}$ we have $\vert\sigma_{s+1}(u)\vert^2=\normover{F}{E}{u}$. This means that the rank of the group
\[
\{\vert\sigma_{s+1}(u)\vert^2 \st u\in\unitpos{F}\}
\]
is the rank of the image of the norm map \map{\text{Nm}_{F/E}}{\unitpos{F}}{\unitpos{E}}. But $(\unitpos{E})^2\subset\Im(\text{Nm}_{F/E})$, thus $\text{Nm}_{F/E}$ has the same rank as $\unitpos{E}$, which is $n/2+0-1=n/2-1$ by Dirichlet Unit Theorem.
\end{D}

\begin{re}
Theorem~\ref{propodd} holds for the general Oeljeklaus-Toma manifold $M(F,U)$, whenever $U\subset\unit{F}$ has finite index.
\end{re}

\begin{re}
Theorem~\ref{propodd} shows that \cite[Example~2.13]{GOPRVS} holds only for some Oeljeklaus-Toma manifolds.
\end{re}

The following two examples describe the case $[F:\Q]$ even.
\begin{ex}\label{maximal}
Pick monic polynomials $f_1$, $f_2$ and $f_3$ in $\Z[X]$ of degree $2n$ such that:
\begin{itemize}
\item $f_1$ is irreducible modulo $2$;
\item $f_2$ splits as a product of a linear factor and an irreducible polynomial modulo $3$;
\item $f_3$ is a product of an irreducible polynomial of degree $2$ and of two irreducible polynomials of odd degree modulo $5$.
\end{itemize}
Then for every polynomial $g\in \Z[X]$ of degree $2n-1$ the polynomial
\[
f=-15f_1+10f_2+6f_3+30g
\]
is monic, is irreducible (since its reduction modulo $2$ is irreducible), and has maximal Galois group $S_{2n}$ (proceed as in \cite[Example 4.31]{MilFGT}, noting that $30\equiv 0$ modulo $2$, $3$ and $5$). For suitable choices of $g$ we will have that $f$ has exactly $2$ non-real roots (proceed as in \cite[Remark~1.1]{OeTNKC}, observing that the set of polynomials $\{-15f_1+10f_2+6f_3+30g, \deg g=2n-1\}$ is a lattice in $\Q^{2n}$). Let $F_f$ be the splitting field of $f$ and fix an isomorphism between $\Gal(F_f/\Q)$ and $S_{2n}$; let $H\subset\Gal(F_f/\Q)$ be the subgroup corresponding to $S_{2n-1}$ viewed as the set of all permutations fixing $1$. Then $F_f^{H}$ has no proper subfields as $S_{2n-1}\subset S_{2n}$ is a maximal subgroup, and by Theorem~\ref{propodd}, point~\ref{even}, the corresponding Oeljeklaus-Toma manifold $M(F_f^{H},\unitpos{F_f^{H}})$ has maximal rank $\lckrk{M}=b_1(M)=2n-2$.
\end{ex}

\begin{ex}\label{half}
Pick an arbitrary totally real number field $E$ of degree $n$. Let $\alpha$ be a primitive element of $E$ over $\Q$ and let $\alpha_1=\alpha,\alpha_2,\dots,\alpha_n$ be the conjugates of $\alpha$: we can assume $\alpha_1>\alpha_2>\dots>\alpha_{n}$. Let $\sigma_i$ be the embedding corresponding to $\alpha_i$, and let $q\in\Q$ such that $\alpha_{n-1}>q>\alpha_n$. Take $F\ug E(\sqrt{\alpha-q})$. Then $[F:E]=2$ (otherwise, $\alpha-q=e^2$ for some $e\in E$: but then $\sigma_n(\alpha)-q=\sigma_n(e^2)$ so $\alpha_n-q>0$ since $\sigma_n(e)\in\R$ as $E$ is totally real), and $F$ admits exactly 2 complex non-real embeddings (the $[F:E]$ extensions of $\sigma_n$ to $F$). Then by Theorem~\ref{propodd}, point~\ref{even}, the corresponding Oeljeklaus-Toma manifold $M(F,\unitpos{F})$ has rank $\lckrk{M}=\frac{b_1(M)}{2}=n-1$.
\end{ex}

\begin{re}
Example~\ref{half} relies on the existence of a totally real number field of degree $n$, for an arbitrary $n$. This can be shown this way. First, recall that if $p$ is an arbitrary prime number, and $\zeta$ is be a primitive root of unity of order $p$, then $\Q\subset\Q(\zeta)$ is a Galois extension of degree $p-1$, with Galois group cyclic of order $p-1$.
Moreover, $\Q\subset\Q(\zeta+\frac{1}{\zeta})$ is a totally real Galois extension of $\Q$, with Galois group cyclic of order $\frac{p-1}{2}$. Now, choose a prime $p$ such that $n$ divides $\frac{p-1}{2}$ (Dirichlet's Theorem on prime numbers in arithmetic progressions), and choose a subgroup $H$ of $\Gal(\Q(\zeta+\frac{1}{\zeta})/\Q)$ of index $n$. Then $\Q(\zeta+\frac{1}{\zeta})^H$ is a subfield of $\Q(\zeta+\frac{1}{\zeta})$, hence totally real, and $\Q\subset\Q(\zeta+\frac{1}{\zeta})^H$ is Galois.
Thus $[\Q(\zeta+\frac{1}{\zeta})^H:\Q]$ is the cardinality of its Galois group, that is exactly $n$.
\end{re}

\begin{re}
We summarize in the following table what we presently know about relations between $\lckrk{M}$ and properties of $M$. If $M=(K,\Gamma)$, by ``Potential'' in the Table we mean there exists a global K\"ahler potential on $K$, and by ``Automorphic potential'' we mean there exists a global K\"ahler potential on $K$ such that $\Gamma$ acts on it by homotheties (see \cite{OrVMNC,OrVLCK}).
\end{re}

\begin{small}
\begin{center}
\begin{tabular}{|m{.4\textwidth}|c|m{.3\textwidth}|}
\hline
\centering \textbf{Statement} & \textbf{True/False} & \textbf{Proof or Refutation} \\ \hline
\centering $b_1=1\Rightarrow\lckrk{M}=1$ & True & $1\le\lckrk{M}\le b_1$\\ \hline
\centering $\lckrk{M}=1\Rightarrow b_1=1$ & False & Induced Hopf bundles over curves of large genus in $\CP{2}$: $\lckrk{M}=1$, arbitrarily large $b_1$\\ \hline
\centering Vaisman $\Rightarrow\lckrk{M}=1$ & True & \cite[Corollary~4.7]{GOPRVS} \\ \hline
\centering Automorphic potential $\Rightarrow$ $M$ can be deformed to $M'$, with $\lckrk{M'}=1$ & True & \cite[Proposition~2.5]{OrVLCK} and \cite[Proofs of 5.2 and 5.3]{OrVMNC} \\ \hline
\centering $\lckrk{M}=1\Rightarrow$ Automorphic potential & False & Inoue surfaces \\ \hline
\centering Potential $\Rightarrow\lckrk{M}=1$ & False & Oeljeklaus-Toma manifolds as in Theorem \ref{propodd}, with $s>2$: $\phi$ is a potential, and $\lckrk{M}=b_1$ or $b_1/2$ \\ \hline
\centering $\lckrk{M}=1\Rightarrow$ Potential & False & Diagonal Hopf surface blown up in one point: $b_1=1$, so $\lckrk{M}=1$, and no potential because the universal covering contains compact complex submanifolds\\ \hline
\centering $\exists M$ with $\lckrk{M}=b_1>1$ & True & Take $n=2$, $f_1=X^4+X+1=f_2$, $f_3=(X^2+2)(X-1)(X+1)$ and $g=0$ in Example \ref{maximal} \\ \hline
\centering $\exists M$ with $b_1>1$ and $\lckrk{M}=b_1/2$ & True & Take $n=2$, $L=\Q(\sqrt{2})$, $\alpha=\sqrt{2}$ and $q=0$ in Example \ref{half} \\ \hline
\end{tabular}
\end{center}
\end{small}

\section{Algebraic Number Theory background}\label{background}

Let $F$ be a number field, that is, a finite extension of $\Q$. Such an extension is algebraic \cite[Proposition 1.30]{MilFGT}, that is, elements $x$ in $F$ satisfy $P(x)=0$, where $P$ is a polynomial in $\Q[X]$. Whenever $P$ can be chosen monic and with coefficients in $\Z$, $x$ is said to be an \emph{algebraic integer}, and all algebraic integers in $F$ form a ring usually denoted by \integer{F}:
\[
\integer{F}\ug\{x\in F \st x^l+a_{l-1}x^{l-1}+\dots+a_1x+a_0=0,\quad a_i\in\Z\}
\]

Algebraic integers $\integer{F}$ are a Dedekind domain \cite[Theorem 3.29]{MilANT}. Moreover, it is well-known that $\integer{F}/I$ is a finite ring whenever $I$ is a proper ideal in $\integer{F}$.

For any $x\in F$, there is one and only one monic, irreducible $P_x\in\Q[X]$ such that $P_x(x)=0$. Such a $P_x$ is called the \emph{minimal polynomial} of $x$ over $\Q$, and algebraic integers are characterized by having minimal polynomial in $\Z[X]$ \cite[Proposition 2.11]{MilANT}:
\[
\integer{F}=\{x\in F \st P_x\in\Z[X]\}
\]
The quotient ring $\Q(x)=\Q[X]/(P_x)$ is a field because $P_x$ is irreducible, and it is the smallest field containing $\Q$ and $x$. By \cite[Primitive Element Theorem, 5.1]{MilFGT}, every number field is obtained this way, so it is not restrictive to think of $F$ as $\Q(\alpha)$, for a fixed $\alpha\in\C$. The degree $\deg P_\alpha=n=[F:\Q]$ of $P_\alpha$ is then the dimension of $F$ as a vector space over $\Q$, a basis for $F$ being $\{1,\alpha,\dots,\alpha^{n-1}\}$. Any intermediate field between $F$ and $\Q$ is of the form $\Q(x)$, for some $x\in F$, and $\deg P_x$ is a divisor of $\deg P_\alpha$ \cite[Proposition~1.20]{MilANT}. The ring of algebraic integers $\integer{F}$ is a free $\Z$-module of rank $n$ \cite[Page 29]{MilANT}.

Any root $\alpha_i$ of $P_\alpha$ induces a field embedding \map{\sigma_i}{F}{\C}, by
\[
x_0+x_1\alpha+\dots+x_{n-1}\alpha^{n-1}\longmapsto x_0+x_1\alpha_i+\dots+x_{n-1}\alpha_i^{n-1}
\]
Clearly, $\sigma_i(x)=x$ for every $x\in\Q$, and the $\sigma_i$ are the only embeddings of $F$ into $\C$ with this property \cite[Proposition 2.1b]{MilFGT}.
Moreover, $\sigma_i(F)\subset\R$ if and only if $\alpha_i\in\R$, and $F$ is called \emph{totally real} if $\sigma_i(F)\subset\R$ for all $i$.

If $E$ is any intermediate field $F\supset E\supset\Q$, we can consider the finite group of all automorphisms of $F$ fixing $E$ pointwise, denoted by $\Aut{F/E}$. Then, for any subgroup $H$ of $\Aut{F/E}$, we have the subfield of $F$ given by $F^H=\{x\in F\st Hx=x\}$. The key point of Galois Theory is that this ``subfield-subgroup-subfield'' correspondence is $1-1$, for a certain class of extensions $F$, called \emph{Galois extensions} \cite[Theorem~3.16]{MilFGT}.

Galois extensions are characterized by the following equivalent conditions \cite[Theorem~3.10]{MilFGT}:
\begin{itemize}
\item $F=E(\alpha)$ contains all roots of $P_\alpha$, where $P_\alpha$ is the minimal polynomial of $\alpha$ \emph{over $E$}.
\item $F$ contains all roots of an irreducible polynomial $P\in E[X]$ (we say that $F$ is the \emph{splitting field} of $P$).
\item $\Aut{F/E}$ contains $n$ elements, where $n=\deg P_\alpha=[F:E]$.
\item $E=F^{\Aut{F/E}}$.
\end{itemize}

Whenever $F$ is a Galois extension of $E$, the group $\Aut{F/E}$ is called \emph{Galois group of $F$ over $E$}, and it is denoted by $\Gal(F/E)$. If $x\in F$ and $g\in\Gal(F/E)$, then $g(x)$ is called a \emph{Galois conjugate of $x$} (briefly, a conjugate of $x$). One of the many nice properties of Galois extensions is that for any $x\in F$ one has $e\ug\prod_{g\in\Gal(F/E)}g(x)\in E$: this appears evident from the fact that $e$ is fixed by $\Gal(F/E)$ and the fact that $E=F^{\Gal(F/E)}$.

If $F=\Q(\alpha)$ and \map{\sigma_i}{F}{\C} are defined as above, we see that $\Q\subset F$ is Galois if and only if $\alpha_i\in F$, and in this case $\sigma_i\in\Gal(F/\Q)$. This implies that $\Q\subset F$ in Section~\ref{secpot} is never a Galois extension, since there are both real and complex non-real embeddings.

An example of Galois extension is $E\subset F$ in proof of Theorem~\ref{propodd}, point~\ref{even}, since $[F:E]=2$. The non-trivial element of $\Gal(F/E)$ is the complex conjugation.

Multiplication by any $x\in F$ can be viewed as a $\Q$-linear map on $F$: the \emph{norm} of $x\in F$, denoted by $\norm{x}$, is the determinant of this linear map. Since the characteristic polynomial $c_x$ of $x$ as a linear map is related to the minimal polynomial $P_x$ by \cite[Proposition 5.40]{MilFGT}
\begin{equation}\label{char}
c_x=P_x^{[F:\Q(x)]}
\end{equation}
one can prove that, for any $x\in F$, one has \cite[Remark 5.43]{MilFGT}
\begin{equation*}
\norm{x}=\prod_{i=1,\dots,n}\sigma_i(x)
\end{equation*}

The norm can be defined the same way as above for a field extension $E\subset F$. One obtains this way a map \map{\text{Nm}_{F/E}}{F}{E}.

The norm can be used to distinguish elements of the multiplicative group $\unit{F}$ of units of $\integer{F}$, using the following result: $x\in \integer{F}$ is a unit if and only if $\norm{x}=\pm 1$ \cite[Lemma 5.2]{MilANT}.

A \emph{positive unit} is a unit which is positive in all real embeddings of $F$: if $n=s+2t$, with $s$ the number of real roots, and $t$ the number of conjugate pairs of complex roots of $P_\alpha$ (with $\alpha_{s+i}$ the complex conjugate of $\alpha_{s+i+t}$), one defines
\begin{equation}\label{norm}
\unitpos{F}\ug\{u\in\unit{F}\st\sigma_i(u)>0 \text{ for }i=1,\dots,s\}
\end{equation}
Thus, for any positive unit $u$, one has:
\begin{equation}\label{prodsigma}
\prod_{i=1,\dots,n}\sigma_i(x)=1
\end{equation}

It is a classical fact that the norm takes units to units, and positive units to positive units:
\[
\text{Nm}_{F/E}|_{\unit{F}}:\unit{F}\longrightarrow\unit{E},\qquad\text{Nm}_{F/E}|_{\unitpos{F}}:\unitpos{F}\longrightarrow\unitpos{E} 
\]

The units group $\unit{F}$ is a finitely generated abelian group with rank $s+t-1$ \cite[Dirichlet Unit Theorem, 5.1]{MilANT}. Its torsion is the set of roots of $1$ contained in $\integer{F}$, so whenever $s>0$ (that is, $P_\alpha$ has at least one real root), it must be $\{\pm1\}$.
Clearly, $\unitpos{F}$ doesn't contain $-1$, so it is free. Moreover, it contains the subgroup $\{u^2\st u\in\unit{F}\}$, which has finite index in $\unit{F}$. Thus, $\unitpos{F}$ has rank $s+t-1$.

\vspace{1em}

\noindent{\bf Acknowledgements.} The authors thanks Rosa Gini and Liviu Ornea for their valuable support, Matei Toma for important clarifications regarding \cite{OeTNKC}, and Antonio J.~Di Scala for pointing us an erroneous argument used in the proof of Theorem~\ref{mainoe}.


\end{document}